\documentclass[11pt,
]{article}

\usepackage{amssymb,amsfonts,amstext,amsmath,amsthm,
latexsym,mathrsfs,amsbsy}
\usepackage{bbold}

\usepackage{marginnote}

\usepackage{mparhack}

\makeatletter
\if@twoside%
\else\fi%
\makeatother  

\newcommand{\nek}{\newcommand}
\nek{\renek}{\renewcommand}
\nek{\vyk} [1] {}

\nek{\ubf}{\fontseries{b}\selectfont}
\nek{\bfit}{\bfseries\itshape}
\nek{\bftt}{\ttfamily\bfseries\upshape\selectfont}

\nek{\parf}{\subsection}
\renek{\thesubsection}{\arabic{subsection}}
\nek{\punk}{\subsubsection}
\renek{\thesubsubsection}
{\thesubsection\alph{subsubsection}}

\theoremstyle{plain}

\newtheorem{theorem}             {Theorem} [subsection]
\newtheorem{corollary}  [theorem]{Corollary}
\newtheorem{prop}  [theorem]{Proposition}
\newtheorem{lemma}      [theorem]{Lemma}
\newtheorem{slemma}      [theorem]{Sublemma}
\newtheorem{claim}      [theorem]{Claim}

\theoremstyle{definition}

\newtheorem{definition} [theorem]{Definition}
\newtheorem{vopr}       {Problem}
\newtheorem{prim}       [theorem]{Example}
\newtheorem{zam}        [theorem]{Remark}
\newtheorem*{prF}   {Proof}    
\newtheorem*{ack}   {Acknowledgement}    
\newtheorem{ggi}   [theorem] {Соглашение}

\nek{\thsp}{\hspace{0.1ex}}
\nek{\bvo}{\begin{vopr}}
\nek{\evo}{\end{vopr}}
\nek{\back}{\begin{ack}}
\nek{\eack}{\end{ack}}
\nek{\bpro}{\begin{prop}}
\nek{\epro}{\end{prop}}
\nek{\bcor}{\begin{corollary}}
\nek{\ecor}{\end{corollary}}
\nek{\bdf} {\begin{definition}}
\nek{\edf} {\qed\end{definition}}
\nek{\eDf} {\end{definition}}
\nek{\bgg} {\begin{ggi}}
\nek{\egg} {\qed\end{ggi}}
\nek{\ble} {\begin{lemma}}
\nek{\ele} {\end{lemma}}
\nek{\bsl} {\begin{slemma}}
\nek{\esl} {\end{slemma}}
\nek{\bcl} {\begin{claim}}
\nek{\ecl} {\end{claim}}
\nek{\bpri}{\begin{prim}}
\nek{\epri}{\qed\end{prim}}
\nek{\eprI}{\end{prim}}
\nek{\bte} {\begin{theorem}}
\nek{\ete} {\end{theorem}}
\nek{\baq} {\begin{aq}}
\nek{\eaq} {\end{aq}}
\nek{\bre}{\begin{zam}}
\nek{\ere}{\qed\end{zam}}
\nek{\bpf} {\begin{prF}} 
\nek{\epf} {\qed\end{prF}} 
\nek{\epg} {\end{prF}} 
\nek{\epF}[1]{\hfill\hbox{$\square$\ ({\small#1})}\end{prF}}

\nek{\qeD}[1]{\hfill\hbox{$\square$\ ({\small#1})}}

\nek{\ben}{\begin{enumerate}\itemsep=0.2em}
\nek{\een}{\end{enumerate}}
\nek{\bde}{\begin{description}\itemsep=0.2em}
\nek{\ede}{\end{description}}
\nek{\bit}{\begin{itemize}\itemsep=0.2em}
\nek{\eit}{\end{itemize}}
\nek{\bay}{\begin{array}}
\nek{\eay}{\end{array}}
\nek{\bce}{\begin{center}}
\nek{\ece}{\end{center}}
\nek{\bqu}{\begin{quotation}\noindent}
\nek{\equ}{\end{quotation}}
\nek{\lh} [1] {\mathop{\tt lh}(#1)}
\nek{\card} {\mathop{\tt card}}
\nek{\dom} {\mathop{\tt dom}}
\nek{\clos} {\mathop{\tt clos}}
\nek{\num} {\mathop{\tt num}}
\nek{\ran} {\mathop{\tt ran}}
\nek{\rank} {\mathop{\tt rank}}
\nek{\tsup} {\mathop{\tt sup}}
\nek{\tmax} {\mathop{\tt max}}
\nek{\sep}{{\tt Sep}}
\nek{\red}{{\tt Red}}
\nek{\TC}  {{\rm TC}\hspace{0.4ex}}
\nek{\hc}  {\mathrm{HC}}
\nek{\rL} {\text{\ubf L}}
\nek{\rS} {{\mathbf S}}
\nek{\rV} {\text{\ubf V}}

\nek{\ZFC} {\text{\ubf ZFC}}
\nek{\zfc} {\ZFC}
\nek{\zfcm} {\zfc^-}
\nek{\al}  {\alpha}
\nek{\ga}  {\gamma}
\nek{\Ga}  {\Gamma}
\nek{\da}  {\delta}
\nek{\Da}  {\Delta}
\nek{\kpa} {\kappa}
\nek{\la}  {\lambda}
\nek{\ve}  {\varepsilon}
\nek{\vpi} {\varphi}
\nek{\vpo} {\overline\vpi}
\nek{\sg}  {\sigma}
\nek{\Sg}  {\Sigma}
\nek{\om}  {\omega}
\nek{\Om}  {\Omega}
\nek{\lom} {^{<\om}} 
\nek{\omi} {\om_1}
\nek{\oli} {\omi^\rL}
\nek{\omb} {\om_2}
\nek{\omil}{\om_1^\rL}
\nek{\za}  {\zeta}
\nek{\ali} {\aleph_1}
\nek{\ald} {\aleph_2}
\nek{\alil} {\aleph_1^\rL}
\nek{\fs}[2]{{\mathbf\Sigma}^{#1}_{#2}}
\nek{\fp}[2]{{\mathbf\Pi}^{#1}_{#2}}
\nek{\fd}[2]{{\mathbf\Delta}^{#1}_{#2}}
\nek{\Gd}{{\mathbf G}_\delta}
\nek{\is}[2]{\varSigma^{#1}_{#2}}
\nek{\ip}[2]{\varPi^{#1}_{#2}}
\nek{\id}[2]{\varDelta^{#1}_{#2}}
\nek{\BBB}{\hspace{0.1ex}}
\nek{\dA}{{\BBB{\mathbb A}\BBB}}
\nek{\dB}{{\BBB{\mathbb B}\BBB}}
\nek{\dK}{{\BBB{\mathbb K}\BBB}}
\nek{\dN}{{\BBB{\mathbb N}\BBB}}
\nek{\dP}{{\BBB{\mathbb P}\BBB}}
\nek{\dQ}{{\BBB{\mathbb Q}\BBB}}

\nek{\uff} {\uar f}
\nek{\uzp} {\uar\zp}
\nek{\uzq} {\uar\zq}
\nek{\udp} {\uar\dP}
\nek{\ump} {\uar\mP}
\nek{\uzd} {\uar\zD}
\nek{\uar} [1] {#1^{\uparrow}}
\nek{\dR}{{\BBB{\mathbb R}\BBB}}
\nek{\dW}{{\BBB{\mathbb W}\BBB}}
\nek{\dX}{{\BBB{\mathbb X}\BBB}}
\nek{\dY}{{\BBB{\mathbb Y}\BBB}}
\nek{\gM}{{\BBB{\goth M}\BBB}}
\nek{\gN}{{\BBB{\goth N}\BBB}}
\nek{\cM}{{\BBB{\skri M}\BBB}}
\nek{\cP}{{\BBB{\skri P}\BBB}}
\nek{\pu}  {\varnothing}
\nek{\sq}  {\subseteq}
\nek{\qs}  {\supseteq}
\nek{\su}  {\subset}
\nek{\sneq}{\subsetneqq}
\nek{\eqv} {\mathbin{\,\Longleftrightarrow\,}}
\nek{\imp} {\mathbin{\,\Longrightarrow\,}}
\nek{\mpi} {\mathbin{\,\Longleftarrow\,}}
\nek{\lra} {\longrightarrow} 
\nek{\we}  {{\mathbin{\hspace*{0.2ex}^\smallfrown}}}
\nek{\sus} {{\exists\,}}
\nek{\kaz} {{\forall\,}}
\nek{\ti}  {\times}
\nek{\dm}  {$$}
\nek{\iy}  {\infty}
\nek{\abs}[1]{|#1|}
\nek{\abt}[1]
{\mathopen{{|}\hspace*{-0.2ex}{|}}#1\mathclose{{|}\hspace*{-0.2ex}{|}}}

\nek{\nin}{\notin}
\nek{\res} {{\hspace*{0.1ex}\restriction\hspace*{0.3ex}}}
\nek{\ares} {{\hspace*{0.1ex}\uparrow\hspace*{0.3ex}}}

\nek{\bang} [1] {\big\langle #1\big\rangle}
\nek{\ang} [1] {\langle #1\rangle}
\nek{\ans} [1] {\{\hspace{0.1ex}#1\hspace{0.1ex}\}}
\nek{\dd}[1]{$\mtho\hspace{0.2ex}{#1}$-\hspace{0.0ex}}
\nek{\itla}{\item\label}
\nek{\itlb} [1] {\itla{#1}\hspace*{0ex}\imar{#1}}%

\nek{\bez}{\smallsetminus}
\nek{\rit} [1] {{\it#1\/}}

\nek{\ens} [2] {\ans{{#1\hspace{0.5ex}{:}}\zz\hspace{0.5ex}#2}}
\nek{\zz} {\linebreak[0]} 
\nek{\goth}{\mathfrak}
\nek{\skri}{\mathscr}

\nek{\yo} {,\linebreak[0]}
\nek{\yi} {\hspace{\mathsurround},\linebreak[0]\hspace{\mathsurround}}
\nek{\yd} {\hspace{\mathsurround},\linebreak[0]\:}
\nek{\yt} {\hspace*{\mathsurround}\text{,}\linebreak[0]\;}

\nek{\viv} {\text{vice versa}}
\nek{\te} {\text{т.\resp е.}}
\nek{\Te} {\text{T.\resp е.}}
\nek{\ie} {\text{\sl i.\resp e.}}
\nek{\eg} {\text{\sl e.\resp g.}}
\nek{\poo} {\text{w.\resp r.\resp t.}}
\nek{\noo} {\text{w.\resp l.\resp o.\resp g.}}

\nek{\iesp}{\hspace{0.3ex}}
\nek{\resp}{\hspace{0.25ex}}

\nek{\itsep}{\itemsep=0.25ex plus 0.1ex minus 0.1ex}

\nek{\tenu}[1]{
\def\theenumi{#1}
\itsep}

\nek{\tenv}[1]{
\def\theenumii{#1}
\itsep}

\nek{\Aenu}{\tenu{(\Alph{enumi})}}
\nek{\aenu}{\tenu{(\alph{enumi})}}
\nek{\aenup}{\tenu{(\alph{enumi}{\mtho$'$})}}
\nek{\nenu}{\tenu{{\rm(\arabic{enumi})}}}
\nek{\nenup}{\tenu{{\rm(\arabic{enumi}{\mtho$'$})}}}

\nek{\fenu}{\itsep\tenu{{\mtho$(\fnsymbol{enumi})$}}}
\nek{\renu}{\itsep\tenu{{\rm(\roman{enumi})}}\itsep}
\nek{\Renu}{\itsep\tenu{{\rm(\Roman{enumi})}}\itsep}

\nek{\atc} {\addtocounter{enumi}1}
\nek{\atm} {\addtocounter{enumi}{-1}}

\nek{\stk} [2] {\ang{#1\hspace{0.3ex};\hspace{0.1ex}#2}}
\nek{\sis} [2] {\ang{#1}_{#2}}
\nek{\sid} [3] {\bang{#1}{}^{#2}_{#3}}

\renek{\cM} {\goth M}

\nek{\lam} [1]
{\label{#1}\hspace*{-3pt}\imar{#1}%
}%
\nek{\las} [1]
{\label{#1}\imae{#1}}%

\nek{\imar}[1]{\marginpar[
\flushright\footnotesize%
$\mtho\longrightarrow$\\%
\vspace{-1ex}{#1}\vspace*{1ex}]%
{
\flushleft\footnotesize%
$\mtho\longleftarrow$\\%
\vspace{-1ex}{#1}\vspace*{1ex}}%
}%

\nek{\imae}[1]{\marginpar[
\flushright\footnotesize\vspace{-4ex}%
$\mtho\longrightarrow$\\%
\vspace{-1ex}{#1}$\mtho$\vspace*{1ex}]%
{
\flushleft\footnotesize\vspace{-4ex}%
$\mtho\longleftarrow$\\%
\vspace{-1ex}{#1}$\mtho$\vspace*{1ex}}
}%

\nek{\pws} [1] {\cP(#1)}
\nek{\lla} {\,\land\,}

\nek{\snos} [1] {\,\footnote{\ #1}}
\nek{\snom}   {\,\footnotemark}
\nek{\snot} [1] {\footnotetext{\ #1}}

\nek{\ap} {{\hspace*{0.2ex}\boldsymbol\cdot\hspace*{0.2ex}}} 

\nek{\onto} {\overset{{\text{\rm onto}}}{\longrightarrow}}

\nek{\La} {\Lambda}
\nek{\alo} {{\aleph_0}}

\nek{\leqv} {\,\eqv\,}

\nek{\dop} [1] {{#1}{}^{\complement}}

\nek{\ba}{\beta}

\nek{\etc} {{\sl etc}.}

\nek{\qand} {\quad\text{and}\quad}

\nek{\vt}{\vartheta}

\nek{\nse} {\om^{<\om}}
\nek{\bse} {2^{<\om}}

\nek{\namx} [1] 
{\overset{\text{\mtho$\hspace*{0.2ex}_\text{\Large\bf.}$}}{#1}}
\nek{\namy} [1] 
{\overset{\text{\mtho$\hspace*{0.5ex}_\text{\Large\bf.}$}}{#1}}
\nek{\namz} [1] 
{\overset{\text{\mtho$\hspace*{0.9ex}_\text{\Large\bf.}$}}{#1}}

\nek{\cD} {\mathscr D}

\nek{\np}{\newpage}

\nek{\sqf} {\sq^{\text{\tt fin}}}
\nek{\sqd} {\sq^{\text{\tt fd}}}

\nek{\bv} {\bigvee}
\nek{\sqfv} {\sq^{\text{\tt fin}}\bigvee}
\nek{\sqdv} {\sq^{\text{\tt  fd}}\bigvee}

\nek{\sqfc} {\sq^{\text{\tt fin}}\bigcup}
\nek{\sqdc} {\sq^{\text{\tt  fd}}\bigcup}

\nek{\ka}{\kappa}

\nek{\pwor} [1] {{\text{\rm PWO}}(#1)}
\nek{\sepa} [1] {{\text{\ubf Sep}}(#1)}
\nek{\redu} [1] {{\text{\ubf Red}}(#1)}

\nek{\pet} {\mathbf{PT}}
\nek{\ret} [1] {\res_{\hspace*{0.05ex}#1}}

\renek{\leq} {\leqslant}
\renek{\geq} {\geqslant}

\nek{\limp} {\mathrel{\,\imp\,}}

\nek{\spe} [1] {\text{\bf SC}(#1)}
\nek{\spf} [1] {\text{\bf SC}(#1)}
\nek{\spg} [1] {\text{\bf SC}^{<\om}(#1)}

\nek{\vys} [1] {\text{\tt hgt}(#1)} 

\nek{\cle} {\preccurlyeq}
\nek{\cls} {\prec}

\renek{\refname} {{\large\bf References}}

\nek{\dphi} {\mathbf\Phi}

\nek{\dU}{\mathbb U}
\nek{\uu}{^{\text{\tt fu}}}

\nek{\Tv}{\overrightarrow{T}}
\nek{\Sv}{\overrightarrow{S}}

\nek{\jLa} {\boldsymbol\La}
\nek{\jta} {\boldsymbol\tau}
\nek{\jsg} {\boldsymbol\sg}
\nek{\jrho} {\boldsymbol\rho}
\nek{\jro} {\jrho}
\nek{\jw} {\boldsymbol w}

\nek{\jv} {\mathbf V}
\nek{\ju} {U}

\nek{\dn} {2^\om}
\nek{\dnn} {{(\dn)}{}^\om}
\nek{\bn} {\om^\om}

\nek{\dox}{{\namy {\boldsymbol x}}}
\nek{\doa}{{\namy {a}}}
\nek{\dob}{{\namy {b}}}
\nek{\dof}{{\namy {f}}}
\nek{\doH}{{\namy {H}}}

\nek{\rc} {\mathbf c}
\nek{\rd} {\mathbf d}
\nek{\rpi} {\dox}

\nek{\kn} [1] {\underline{#1}}

\nek{\dplom} {\dd{\dP\lom}}
\nek{\dplo} {\dd\plo}

\nek{\plo} {\dP\lom}
\nek{\uplo} {(\dP\cup\dU)\lom}

\nek{\cll} {\mathrel{{\cls}\hspace*{-0.9ex}{\cls}}}

\nek{\vpj} [2] {\vpi^{#1}_{#2}}
\nek{\Phj} [2] {\Phi^{#1}_{#2}}

\nek{\vpt} [3] {\vpi^{#1}_{#2#3}}

\nek{\od} {\text{\rm OD}}
\nek{\hod} {\text{\rm HOD}}
\nek{\rod} {\text{\rm ROD}}

\nek{\Eo} {\mathrel{{\text{\sf E}}_0}}

\nek{\zf} {\text{\ubf ZF}}

\nek{\uf} [2] {{U}^\dphi_{#1}(#2)}
\nek{\ufi} [1] {{U}^\dphi_{#1}}

\nek{\tf} [3] {{U}^\dphi_{#1#2}(#3)}
\nek{\tfi} [2] {{U}^\dphi_{#1#2}}

\nek{\tx} [2] {{T}^\dphi_{#1}(#2)}
\nek{\txi} [1] {{T}^\dphi_{#1}}

\nek{\ty} [3] {{T}^\dphi_{#1#2}(#3)}

\nek{\lel} {\leq_{\rL}}

\nek{\kc} [2] {C_{#1#2}}
\nek{\xc} [2] {C'_{#1#2}}
\nek{\kd} [2] {D_{#1#2}}
\nek{\kcp} [3] {C_{#1#2}^{#3}}
\nek{\kr} [3] {\jrho_{#1#2}^{#3}}
\nek{\kR} [2] {R_{#1#2}}

\nek{\jcp} [4] {C_{#1#2}^{#3#4}}
\nek{\jr} [4] {\jrho_{#1#2}^{#3#4}}

\nek{\roo} [1] {\text{\tt stem}(#1)}

\nek{\cI} {\mathcal I}
\nek{\fss} [1] {\text{\bf FSS}(#1)}
\nek{\fsd} [1] {\text{\bf FSS}\lom(#1)}

\nek{\ptf} {\text{\ubf {PF}}}

\vyk{
\nek{\jp} {\mathbb p}
\nek{\bbu} {\mathbb u}
\nek{\bbp} {\mathbb p}
\nek{\bbq} {\mathbb q}
\nek{\bbpu} {\bbp{\lor}\bbu}
}

\nek{\bPi} {\mathbb\Pi}

\nek{\mun} {мультичис}
\nek{\muc} {мультикортеж}
\nek{\mus} {system}
\nek{\Mus} {System}
\nek{\ms} [1] {\text{\ubf Sys}(#1)}
\nek{\mss} [1] {\text{\ubf Sys}}

\nek{\usl} {condition}
\nek{\muf} {multiforcing}
\nek{\Muf} {Multiforcing}
\nek{\mut} {multitree}
\nek{\pep} {perfect product}
\nek{\muk} {multicollage}
\nek{\Mut} {Multitree}
\nek{\mt} [1] {\text{\ubf MT}(#1)}
\nek{\mtt} {\text{\ubf MT}}
\nek{\mtp} {\ensuremath{\text{\ubf MT}^+}}

\nek{\smuf} {submultiforcing}
\nek{\Smuf} {Submultiforcing}

\nek{\et}{\eta}

\nek{\mto} {\mapsto}

\nek{\Ord} {\text{\ubf Ord}}

\nek{\gfz} {W}
\nek{\gfu} {K}

\nek{\fo} {\mathop{\hspace*{0.3ex}\tt forc\hspace*{0.3ex}}}

\nek{\tup}{\textup}

\nek{\mtho}{\mathsurround=0mm}
\nek{\msur}{\hspace*{-1\mathsurround}}
\nek{\dsur}{\hspace{-0.3\mathsurround}}
\nek{\hsur}{\hspace{-0.5\mathsurround}}
\nek{\noi}{\noindent}
\nek{\vom}{\vspace{1mm}}
\nek{\vtm}{\vspace{2mm}}

\nek{\fP} {\bPi}

\nek{\raw} [2] {#1{({\to\hspace*{0.2ex}}#2)}}
\nek{\req} [2] {{#1}\ret{#2}}
\nek{\pes} {\text{\bf{ST}}}
\nek{\pel} {\text{\bf{LT}}}
\nek{\lt}{\pel}
\nek{\nq} [1] {\mathrel{{\sq}_{#1}}}
\nek{\leqs} {\leqslant}
\nek{\mth} [2] {h^{#1}_{#2}}
\nek{\ntt} [3] {T^{#1}_{#2}(#3)}
\nek{\snenu}{\tenu{{\rm(\thesubsection.\arabic{enumi})}}}
\nek{\sct} [2] {\text{\ubf Colg}_{#2}(#1)}
\nek{\sco} [1] {\text{\ubf Colg}(#1)}

\nek{\TS}{\textstyle}

\nek{\zn} [2] {\nu^{#1}_{#2}}
\nek{\zh} [3] {\nu^{#1}_{#2#3}}
\nek{\zt} [3] {\tau^{#1}_{#2#3}}
\nek{\zT} [4] {T^{#1}_{#2#3}(#4)}
\nek{\zd} [3] {T^{#1}_{#2#3}}

\nek{\zc} [2] {{#1}(#2)}

\nek{\pc} [2] {{#1}(#2)}

\nek{\jd} {\mathbf D}

\nek{\zp} {P}
\nek{\zq} {Q}
\nek{\zr} {R}

\nek{\ta} {\tau}

\nek{\ah} {\ibb} 
\nek{\bh} {\dbb} 
\nek{\chh} {\idbb} 

\nek{\dpo} {\rP_{\text{\rm coh}}}
\nek{\mpo} [1] {\mP^{#1}_{\text{\rm coh}}}

\nek{\rD} {\mathrel{\mathsf D}}
\nek{\rav} [1] {\mathrel{\mathsf D}_{#1}}
\nek{\rE} {\mathrel{\mathsf E}}
\nek{\nE} {\mathrel{\not{\rE}}}
\nek{\rF} {\mathrel{\mathsf F}}

\nek{\reb} {\mathrel{\hspace*{0.1ex}\le_{\text{\rm B}}\hspace*{0.1ex}}} 
\nek{\seb} {\mathrel{\hspace*{0.2ex}\le_{\text{\rm B}}\hspace*{0.25ex}}} 

\nek{\unH} {{\underline H}}
\nek{\ua} {{\underline a}}
\nek{\ux} {{\underline x}}
\nek{\uy} {{\underline y}}

\nek{\lap} [1] {``#1''}

\nek{\kw}{\cup^{\text{\tt cw}}}
\nek{\bkw} {\bigcup^{\text{\tt cw}}}

\nek{\app}{\boldsymbol\cdot}
\nek{\ima} [2] {{#1}\text{\hspace*{0.2ex}''}{#2}}
\nek{\imb} [2] {{#1}[{#2}]}
\nek{\HC}  {\hc}
\nek{\eko} [1] {[#1]_{\Eo}}
\nek{\eke} [1] {[#1]_{\rE}}

\nek{\bLa} {\boldsymbol\La}
\renek{\bLa} {\mathbb\La}

\nek{\oi} [1] {\text{\tt spl}(#1)}
\nek{\oin} [2] {\text{\tt spl}_{#2}(#1)}

\nek{\jpi} {\mathbf P}
\nek{\jqo} {\mathbf Q}

\nek{\jpq} {\jpi\kw\jqo}

\nek{\uG} {\underline G}

\nek{\raz} {\mathbb Q}
\nek{\ves} {\mathbb R}
\nek{\nat} {\om} 

\renek{\dA}{A}
\renek{\dB}{B}

\nek{\jf} {\boldsymbol J}

\nek{\rP} {P}
\nek{\rQ} {Q}
\nek{\rR} {R}

\nek{\bP} {\mathbf P}              
\nek{\bQ} {\mathbf Q}                 
\nek{\bR} {\mathbf R}

\nek{\bT} {\mathbf T}                 
\nek{\bS} {\mathbf S}                 
\nek{\bU} {\mathbf U}                 
\nek{\bV} {\mathbf V}                 
\nek{\bW} {\mathbf W}                 
\nek{\bZ} {\mathbf Z}                 

\nek{\sfo} {\lt-forcing}

\nek{\nn} {\mathbb n}

\nek{\qq} [2] {q^{#1}_{#2}}

\nek{\mlto} [1] {\text{\ubf MT}(#1)}
\nek{\mlt} {\text{\ubf MT}}

\nek{\zaf} [2] {\mathbf\Da[#1,#2]}
\nek{\qi} [2] {\boldsymbol\nu_{#1 #2}}
\nek{\bk} [2] {\mathbf k_{#1 #2}}

\nek{\pro} [1] {\phi\obr({#1})}

\nek{\img} [1] {\ima\phi{#1}}

\nek{\qa} [2] {#1\mathbin{%
\hspace*{-1.5pt}{\upharpoonright}\hspace*{-4.5pt}{\downharpoonright}%
}#2}

\renek{\qa} [2] {#1[#2]}

\nek{\eeo} {\mathrel{\equiv_{\Eo}}} 

\nek{\rau} [2] {#1{({\Rightarrow\hspace*{0.2ex}}#2)}}
\nek{\eto} {{j_0}}

\nek{\raj} [3] 
{#1{({\overset{#3}{\longrightarrow}
\hspace*{0.2ex}}#2)}}

\nek{\dre} [2] {\dd{(#1,#2)}регулярн}

\nek{\jq}  {Q}
\nek{\sgi} {\sg'}
\nek{\tai} {\ta'}
\nek{\cp} {U'}

\nek{\prf} {\mathrel{{\sqsubset}_{\text{\tt w}}\hspace*{-0.3ex}}}
\nek{\rf} {\sqsubset}
\nek{\rfa} [1] {\rf_{#1}}
\nek{\rfp} [1] {\rf^{+}_{#1}}
\nek{\dpt} [2] {D_{#1}(#2)}
\nek{\zD} {\text{\ubf D}}
\nek{\abc} [1] {|#1|}
\nek{\wed} {с.\,д.}

\nek{\bssq} {\mathrel{\boldsymbol{\sqsubset}}}
\nek{\brsq} {\mathrel{\boldsymbol\sqsubset^\ast\hspace*{-0.2ex}}}

\nek{\pxf} {\mathrel{{\bssq}^\ast\hspace*{-0.3ex}}}

\nek{\duz} [3] {{#1}^{\abc{#2}}_{#3}}

\nek{\ssa} [1] {\ssq_{#1}}
\nek{\rsa} [1] {\ssq^{\ast}_{#1}}
\nek{\ssq} {\sqsubset}

\nek{\mf} {\text{\ubf MF}}
\nek{\wdi} {weakly disjoint}

\nek{\ssb} [1] {\bssq_{#1}}
\nek{\rsb} [1] {\mathrel{\boldsymbol\sqsubset^\ast_{#1}}}
\nek{\ssc} [1] {\bssq^3_{#1}}
\nek{\sse} [3] {\bssq^{#1}_{#2#3}}

\nek{\kmar}[1]
{\marginnote
[{\scriptsize\rm#1$\Rightarrow$}]%
{{\scriptsize\rm$\Leftarrow$#1}}}%

\nek\NoIndent[1]{%
  \begingroup
  \par
  \parshape0
  #1\par
  \endgroup
}

\nek{\bb} {B}
\nek{\df} {\text{\rm DEF}}

\nek{\ft} [2] {{\mathbf T}^{#1}_{#2}}
\nek{\fl} [1] {{\mathbf L}^{#1}}
\nek{\fls} [3] {\raw{\fl{#1}(#2)}{#3}}

\nek{\jn} {{\boldsymbol\nu}}
\nek{\jm} {{\boldsymbol\mu}}

\nek{\bs} {\mathbf s}
\nek{\bt} {\mathbf t}
\nek{\mc} {\text{\ubf MC}}
\nek{\mn} {\text{\ubf MN}}
\nek{\mlh} [1] {\mathop{\tt mlh}(#1)}

\nek{\vek}{\vv}
\nek{\vjpi}{\vek\jpi}
\nek{\vjqo}{\vek\jqo}
\nek{\vmf} {\vek\mf}
\nek{\vmo} {\vmf_{\omi}}

\nek{\vstf}{\vmf}
\nek{\vsto}{\vmo}

\nek{\gp} {\mathbb P}
\nek{\yp} [1] {\gp(#1)}
\nek{\yyp} [2] {\gp_{#1}(#2)}

\nek{\vdp} {\vv{\gp}}
\nek{\gU}{{\BBB{\goth U}\BBB}}

\nek{\mP} {\text{\mtho\boldmath$\mathfrak S$}}
\nek{\mdp} [1] {\gp_{<#1}}

\nek{\mm} [1] {\cM_{#1}}

\nek{\smf} [1] {\mP(#1)}

\nek{\prol} [1] {\mathrel{\sqsubset_{#1}}}
\nek{\vjR} {\vv{\mathbf R}}

\nek{\cne} [1] {C_{\ne#1}}
\nek{\bc} {\mathbf c}
\nek{\bd} {\mathbf d}
\nek{\jb} {\mathbf b}
\nek{\jba} {{\boldsymbol\ba}}
\nek{\je} {\mathbf e}

\nek{\knf} [1] {\kng_{#1}}
\nek{\kng}  {\text{\rm CCF}}
\nek{\ff} [2] {f^{#1}[#2]}
\nek{\fg} [1] {f^{#1}[\zG]}

\nek{\xx} [1] {\bx[#1]}
\nek{\bx}  {\boldsymbol x}

\nek{\cL} {\mathscr L}
\nek{\lis} {\cL\hspace*{-0,4ex}\is}
\nek{\lip} {\cL\hspace*{-0,3ex}\ip}
\nek{\ccf} {\kng}

\nek{\Fo} [1] {\text{\ubf Forc}(#1)}
\nek{\des} [1] {\text{\ubf Des}(#1)}

\nek{\ib}  {\text{\rm IB}}
\nek{\obr} {^{-1}}

\nek{\ww} {W}

\nek{\ddn} [1] {(\dn){}^{#1}}

\nek{\oz} {\text{ОЗ}}

\nek{\zG} {\mathbb G}

\nek{\gs} [1] {\mP_{#1}}

\nek{\PP} {\mathbf{PP}}
\nek{\eqr}{equivalence relation}
\nek{\cE} {\mathscr E}

\nek{\xe} {\boldsymbol e}
\nek{\xf} {\boldsymbol f}

\renek{\imar} [1] {}
\renek{\imae} [1] {}
\renek{\kmar} [1] {}

\begin{document}

\title
{Canonization of smooth equivalence relations on 
infinite-dimensional perfect cubes\thanks
{Vladimir Kanovei's work was supported in part 
by RFBR grant 17-01-00705. 
Vassily Lyubetsky's work was supported in part 
by RSF grant 14-50-00150.}}

\author 
{Vladimir Kanovei\thanks{IITP RAS and RTU (MIIT),
 \ {\tt kanovei@googlemail.com} --- 
contact author. 
}  
\and
Vassily Lyubetsky\thanks{IITP RAS and Moscow University,
\ {\tt lyubetsk@iitp.ru} 
}}

\date 
{\today}

\maketitle


\begin{abstract}
A canonization scheme for 
smooth equivalence relations on $\dR^\om$ modulo 
restriction to infinite perfect products is proposed. 
It shows that given a pair of Borel smooth 
equivalence relations $\rE,\rF$ on $\dR^\om$, there 
is an infinite perfect product $P\sq\dR^\om$ such 
that either ${\rF}\sq{\rE}$ on $P$, or, for some $j<\om$, 
the following is true for all $x,y\in P$: 
$x\rE y$ implies $x(j)=y(j)$, and 
$x\res{(\om\bez\ans j)}=y\res{(\om\bez\ans j)}$ implies 
$x\rF y$.
\end{abstract}

\vyk{
{\scriptsize
\def\contentsname{\normalsize Contents}
\tableofcontents
}
}

\parf{Introduction}
\las{ko}

The canonization problem can be broadly formulated 
as follows. 
Given a class $\cE$ of mathematical structures $E$, 
and a collection $\cP$ of sets $P$ considered as 
{\it large\/}, or {\it essential\/}, find a smaller 
and better structured subcollection $\cE'\sq\cE$ 
such that for any structure $E\in\cE$ with the 
domain $P$ there is a smaller set $P'\in\cP$, 
$P'\sq P$, such that the restricted substructure 
$E\res P'$ belongs to $\cE'$. 
For instance, the theorem saying that every Borel 
real map is either a bijection or a constant on 
a perfect set, can be viewed as a canonization 
theorem, with $\cE=\ans{\text{Borel maps}}$, 
$\cE'=\ans{\text{bijections and constants}}$,
$\cP=\ans{\text{perfect sets}}$. 

We refer to \cite{ksz} as the background of the 
general canonization problem for Borel and 
analytic \eqr s in descriptive set theory.

Among other results, it is established in 
\cite[Section 9.3]{ksz} (theorems 9.26 and 9.27) 
that if $\rE$ belongs to one of two large families 
of analytic \eqr s\snos
{The first family consists of \eqr s classifiable 
by countable structures, the second of those Borel 
reducible to an analytic P-ideal.}  
on $\dnn$
then there is and infinite perfect product 
$P\sq\dnn$ such that ${\rE}\res P$ is 
{\it smooth\/}, that is, simply there exists a Borel 
map $f:P\to\dn$ satisfying ${x\rE y}\eqv f(x)=f(y)$ 
for all $x,y\in P$. 
The canonization problem for smooth \eqr s themselves 
was not considered in \cite{ksz}.\snos
{{\it We avoid any attempt at organizing the very 
complicated class of smooth equivalence relations\/}, 
\cite[page 232]{ksz}.}  
Theorem~\ref{mt}, the main result of this note, 
contributes to this problem.

\parf{Perfect products}
\las{mud}

We consider sets in $\dnn$. 
Let a \rit{\pep} be any set $P\sq\dnn$, such that 
$P=\prod_{k<\om}P(k)$, where 
$P(k)=\ens{x(k)}{x\in P}$ is the \rit{projection} on 
the \dd kth coordinate, ant it is required that each 
set $P(k)$ is a perfect subset of $\dn$.
Let $\PP$ be the set of all \pep s.
\vyk{
Among Borel \eqr s on $\dnn$, we distinguish 
\rit{multiequalities}, \ie, those of the form $\rav U$, 
where $U\sq\om$, defined as follows: 
$x\rD_U y$ iff $x\res U=y\res U$, \ie, 
$x(k)=y(k)$ for all $k\in U$. 
In particular:\vom 

$x\rD_\om y$ iff $x=y$;\vom 

$x\rD_\pu y$ holds for all $x,y\in\dnn$;\vom 

$x\rav{\ans j} y$ iff $x(j)=y(j)$;\vom  

$x\rD_{\om\bez\ans j} y$ 
iff $x(k)=y(k)$ for all $k\ne j$.\vom  

\noi
}%
To set up a convenient notation, say that 
an \eqr\ $\rE$ on $\dnn$:
\bde
\item[\ubf captures $j\in\om$ on $P\in\PP$:]
if 
$x\rE y$ implies $x(j)=y(j)$ for all $x,y\in P$; 

\item[\ubf is reduced to\/ $U\sq\om$ on\/ $P\in\PP$:]
if $x\res U= y\res U$ implies $x\rE y$ for all $x,y\in P$. 
\ede

\bte
\lam{mt}
If\/ $\rE,\rF$ are smooth 
Borel \eqr s on\/ $\dnn$ then there is 
a \pep\/ $P$ such that 
either\/ ${\rF}\sq{\rE}$ on\/ $P$, or, 
for some\/ $j<\om$, 
$\rE$ captures\/ $j$  on\/ $P$ and\/ 
$\rF$ is reduced to\/ $\om\bez\ans j$ on\/ $P$.
\ete

The two options of the theorem are incompatible on \pep s.

The result can be compared to canonization results 
related to \rit{finite} \pep s and \eqr s defined on 
spaces of the form ${(\dn)}{}^m\yt m<\om$. 
Theorem 9.3 in \cite[Section 9.1]{ksz}
implies that every analytic \eqr\ on ${(\dn)}^m$
coincides with one of the multiequalities $\rav U$, 
$U\sq\ans{0,1,\dots,m-1}$, 
on some perfect product $P\sq {(\dn)}^m$, 
where $x\rav U y$ iff $x\res U=y\res U$.
One may ask whether such a result holds for 
\eqr s on $\dnn$ and accordingly infinite \pep s. 
This answers in the negative, even for smooth 
equivalences. 

\bpri
\lam{22}
Let $\rE$ be defined on $\dnn$ so that $x\rE y$ 
iff $x(0)=y(0)$, and also $x(j+1)=y(j+1)$ for all 
numbers $j$ such that $x(0)(j)=0$. 
That $\rE$ is smooth can be witnessed by the map 
sending each $x\in\dnn$ to $a=f(x)\in\dnn$ defined 
so that $a(k)=x(k)$ whenever $k=0$ or $k=j+1$ and 
$x(0)(j)=0$, and $a(k)(n)=0$ for all other $k$ and 
all $n<\om$. 
That $\rE$ is not equal 
(and even not Borel bi-reducible)
to any $\rav U$ on any \pep\ $P\sq\dnn$ is easy.
\epri

The proof of Theorem~\ref{mt} is based on 
splitting/fusion technique 
known in the theory of iterations and 
products of the 
perfect-set forcing (see, \eg, \cite{balav,nge}).

\parf{Splitting}
\las{pss}

The \rit{simple splitting\/} of a perfect 
set $X\sq\dn$ consists of subsets
$\raw Xi=\ens{x\in X}{x(n)=i}$, $i=0,1$, 
where $n=\lh s$ (the lenght of a string $r\in\bse$), 
and $s=\roo X$ 
is the largest string in $\bse$ satisfying $r\su x$ 
for all $x\in X$. 
Then $X=\raw X0\cup\raw X1$ is a disjoint partition of 
a perfect set $X\sq\dn$ onto two perfect subsets.
Splittings can be iterated. 
We let $\raw X\La=X$ for the empty string $\La$, 
and if
$s\in 2^n$, $s\ne\La$ then we define   
$$
\raw Xs=\raw{\raw{\raw{\raw X{s(0)}}{s(1)}}{s(2)}\dots}{s(n-1)}
\,.
$$
If $X,Y\sq\dn$ are perfect sets and $n<\om$ then define 
$X\nq n Y$
(\rit{refinement}), 
if 
$\raw Xs\sq\raw Ys$ for all $s\in2^n$;
$X\nq 0 Y$ is equivalent to $X\sq Y$.
Clearly $X\nq{n+1} Y$ implies $X\nq n Y$
(and $X\sq Y$).

\ble
\lam{nadd}
If\/ $X\sq\dn$ is a perfect set, $s_0\in 2^n$, and\/ 
$A\sq \raw X{s_0}$ is a perfect set,  
then\/ $Y=A\cup\bigcup_{u\in2^n,u\ne s}\raw X{u}$  
is perfect, 
$Y\nq{n} X$, and\/ $\raw{Y}{s_0}=A$.\qed 
\ele       

Now we extend the splitting technique to \pep s. 

\bdf
\lam{phi} 
Fix once and for all 
a function $\phi:\om\onto\om$ taking each 
value infinitely many times, so that 
if $j<\om$ then the following set is infinite:
$$
\pro j=\ens{k}{\phi(k)=j}=
\ans{\bk0j<\bk1j<\bk2j<\ldots<\bk lj<\ldots}.
$$
If $m<\om$ then let $\qi mj$ be   
\kmar{qi mj}%
the number of indices $k<m$, $k\in\pro j$. 
\edf

Let $m<\om$ and $\sg\in2^m$ (a stringh of length $m$). 
If $j\in\img m=\ens{\phi(k)}{k<m}$, 
then the set $\pro j$ 
\kmar{pro n}%
cuts in $\sg$ a substring $\qa\sg j\in2^{\qi mj}$, of 
\kmar{qa sg j}%
length $\lh{\qa\sg j}=\qi mj$,  
defined by $(\qa\sg j)(\ell)=\sg(\bk \ell j)$ for all 
$\ell<\qi mj$. 
Thus the string $\sg\in2^m$ splits in an array   
of strings $\qa\sg j\in2^{\qi mj}$ 
($j\in\img m$) 
of total length $\sum_{j\in\img m}\qi mj=m$.

Let $P$ is a \pep. 
If $j<\om$, $i=0,1$ then 
define a \pep\ $Q=\raj P ij$ so that 
$Q(k)=P(k)$ for all $k\ne j$, but 
$q(j)=\raw{P(j)}i$.
If $\sg\in2^m$ then define a \pep\ 
$\rau P\sg$ by induction so that 
$\rau P\La=P$ and 
$\rau P{\sg\we i}=\raj{\rau P\sg}i\eto$,
where $\eto=\phi(m)\yt m=\lh\sg$.
Note that 
${\rau\zp\sg}(j)=\raw{\zp(j)}{\qa\sg j}$
for all $j$.
In particular, if $j\nin \img m$
then
${\rau\zp\sg}(j)=\zp(j)$, 
because $\lh{\qa\sg j}=\qi mj=0$ holds 
provided $j\nin\img m$.


Let $\zp,\zq$ be \pep s. 
Define $\zp\nq m\zq$, if  
$\zp(j)\nq{\qi mj}\zq(j)$ for all $j$. 
This is equivalent to 
$\rau\zp\sg\sq\rau\zq\sg$ for all $\sg\in2^m$. 

If $\sg,\ta\in 2^m$ then let  
$\,
\zaf\sg\ta=
\om\bez\ens{\phi(i)}{i<m\land \sg(i)\ne\ta(i)} 
$.  

\ble
\lam{nadm}
Let\/ $P\sq\dnn$ is a \pep\ and\/ $m<\om$. Then$:$
\ben
\renu
\itla{nadm1}
if\/ $\sg,\ta\in2^m$, then\/ 
$\rau P\sg\res\zaf\sg\ta=\rau P\ta\res\zaf\sg\ta\;;$

\itla{nadm2}  
if\/ $\sg_0\in 2^m$ and\/  
$B\sq\rau P{\sg_0}$ is a \pep, 
then there is a \pep\/ $Q\nq m P$ satisfying\/
$\rau{Q}{\sg_0}=B\;;$  

\itla{nadm3}  
if\/ $B$ is clopen in\/ $P$ in \ref{nadm2} 
then such a\/ $Q$ 
can be chosen to be clopen in\/ $P\,;$  

\itla{nadm4}  
if\/ $\sg_0,\ta_0\in 2^m$,  
$B\sq\rau P{\sg_0}$ and\/ $B'\sq\rau P{\ta_0}$ 
are \pep s, and\/ 
$B\res\zaf\sg\ta=B'\res\zaf\sg\ta$,
then there is a \pep\/ $R\nq m P$ satisfying\/
$\rau{R}{\sg_0}=B$ and\/ $\rau{R}{\ta_0}=B'$.
\een
\ele
\bpf
\ref{nadm2} 
Apply Lemma~\ref{nadd} componentwise with 
$A=B(j)$ for each $j<\om$.
Namely if $j=\phi(k)$, $k<m$, $\nu=\qi mj$, 
$s_0=\qa{\sg_0} j\in 2^\nu$, then we put 
$Q(j)=
B(j)\cup\bigcup_{s\in2^\nu,s\ne s_0}\raw{P(j)}s$, 
while if $j\nin\img m$ then simply $Q(j)=B(j)$.

\ref{nadm4} 
We first apply \ref{nadm2} for $B\sq\rau P{\sg_0}$, 
getting a \pep\ $Q\nq m P$ such that $\rau Q{\sg_0}=B$. 
We claim that $B'\sq \rau Q{\ta_0}$. 
Indeed if $j\in\zaf{\sg_0}{\ta_0}$ then 
still $\qa{\ta_0} j=s_0$, hence 
$\rau Q{\ta_0}(j)=\raw{Q(j)}{s_0}=B(j)$
by construction, 
therefore $B'(j)=B(j)=\rau Q{\ta_0}(j)$. 
If $j\nin\zaf{\sg_0}{\ta_0}$ then the string 
$t_0=\qa{\ta_0} j\in 2^\nu$ differs from $s_0$, hence 
$\rau Q{\ta_0}(j)=\raw{Q(j)}{t_0}=\raw{P(j)}{t_0}$
by construction, 
therefore $B'(j)\sq \raw{P(j)}{t_0}=\rau Q{\ta_0}(j)$ 
anyway.  
Thus indeed $B'\sq \rau Q{\ta_0}$.

Now we apply \ref{nadm2} for $B'\sq \rau Q{\ta_0}$, 
getting a \pep\ $R\nq m Q$ such that $\rau R{\ta_0}=B'$.
And $\rau R{\sg_0}=B$ holds by the same 
reasons as above.
\epf

\parf{Fusion}
\las{pff}

We begin with a basic fusion lemma, rather 
elementary.

\ble
[fusion]
\lam{fus}
Let\/
$\dots \nq 4 X_3\nq 3 X_2\nq 2 X_1\nq 1 X_0$ be an 
infinite sequence of perfect sets\/ $X_n\sq\dn$.
Then\/  
$X=\bigcap_nX_n$ is perfect and\/ 
$X\nq{n+1}X_n$, $\kaz n$.\qed
\ele
\vyk{
\bpf
Note that  
$\oi T=\ens{\oin{T_n}{n}}{n<\om}$; 
this implies both claims.
\epf
}

A version for \pep s follows:

\ble
[applying Lemma \ref{fus} componentwise]
\lam{fusm}
Let
$\dots\nq5\zp_4\nq4\zp_3\nq3\zp_2\nq2\zp_1\nq1\zp_0$ 
be a sequence of \pep s.
Then\/ $Q=\bigcap_n\zp_n$ is a \pep, 
$Q(j)=\bigcap_m\zp_m(j)$ for all\/ $j<\om$, 
and\/ $Q\nq{m+1}\zp_m$ 
for all\/ $m$.\qed
\ele


\bcor
[see Proposition 9.31 in {\cite[Section 9.3]{ksz}}]
\lam{bos}
If\/ $P\sq\dnn$ is a \pep\ and\/ $B\sq P$ a Borel set 
then there is a \pep\ $Q\sq P$ such that\/ 
$Q\sq B$ or\/ $Q\cap B=\pu$.\qed
\ecor
\vyk{
\bpf
[= Proposition 9.31 in {\cite[Section 9.3]{ksz}}]
$B$ has the Baire property inside $P$, and hence there 
is a string $\sg\in\bse$ such that the \pep\ 
$P'=\rau P{\sg_0}$ satisfies: either $B\cap P'$ is 
meager in $P'$ or $P'\bez B$ is meager in $P'$.
We consider the first case, the other one is pretty 
similar. 
The goal is to define a \pep\ $Q\sq P'$ such that\/ 
$Q\cap B=\pu$.
As $B\cap P'$ is meager, we can assume that 
$B\sq P'\bez\bigcup_mG_m$, where each $G_m\sq P'$ 
is a set relatively open dense in $P'$. 

{\ubf Fact 1.} 
If $R\in\PP$, $R\sq P'$ is clopen 
in $P'$, $\sg\in\bse$, and $m<\om$, 
then there is a string $\tau\in\bse$ 
such that $\sg\sq\tau$ and $\rau{R}\tau\sq G_m$. 
Indeed $\rau{R}\sg$ is clopen in $R$ and in $P'$,
thus $\rau{R}\sg\cap G_m\ne \pu$. 
Let $x\in \rau{R}\sg\cap G_m$. 
Then 
$\ans{x}=\bigcap_{k\ge \lh{\sg_0}}\rau{R}{a\res k}$, 
where $\sg_0\su a\in\dn$. 
Put $\tau=a\res k$ for $k$ large enough.

{\ubf Fact 2.} 
If $R\in\PP$, $R\sq P'$ is clopen 
in $P'$, and $m<\om$, 
then there is a \pep\ $Q\sq R\cap G_m$ 
still clopen in $R$, such that $Q\nq m R$.  
Indeed take some $\sg_0\in2^m.$ 
By Fact 1 there is a \pep\ 
$Z=\rau R{\tau}\sq \rau R{\sg_0}\cap G_m$ clopen in $R$.
By Lemma~\ref{nadm} there is a \pep\/ $R'\nq m R$ with\/
$\rau{R'}{\sg_0}=Z\sq G_m$, still clopen in $R$.
Take another $\sg_1\in 2^m$ and get a \pep\ 
$R''\nq m R'$ with\/
$\rau{R''}{\sg_1}\sq G_m$, still clopen in $R'$ and 
hence in $R$.
Note that $\rau{R''}{\sg_0}\sq\rau{R'}{\sg_0}$ since 
$R''\nq m R'$, thus $\rau{R''}{\sg_0}\sq G_m$ is 
preserved. 
Take another $\sg_2\in 2^m,$ and so on. 
In the end we get a \pep\/ $Q\nq m R$ with\/
$\rau{Q}{\sg}\sq G_m$ for all $\sg\in 2^m,$
hence $Q\sq G_n$, as required.

{\ubf Final construction.} 
Using Fact 2, define a sequence 
$\dots\nq3\zp_2\nq2\zp_1\nq1\zp_0\sq P$ 
as in Lemma~\ref{fusm}, such that $P_m\sq G_m$ 
for all $m$. 
Then the limit \pep\ $Q=\bigcap_mP_m$ satisfies 
$Q\cap B=\pu$, as required.
\epf
}%

\bcor
\lam{bof}
If\/ $P\sq\dnn$ is a \pep\ and\/ $f:P\to\dn$ a Borel map 
then there is a \pep\ $Q\sq P$ such that\/ 
$f\res Q$ is continuous.
\ecor
\bpf
If $n<\om$ and $i=0,1$ then let 
$B_{ni}=\ens{x\in P}{f(x)(n)=i}$. 
Using Corollary~\ref{bos} and Lemma~\ref{nadm}, 
we get a sequence 
$\dots\nq3\zp_2\nq2\zp_1\nq1\zp_0\sq P$ of \pep s
as in Lemma~\ref{fusm}, such that if $m<\om$ and 
$\sg\in 2^m$ then $\rau {P_m}\sg\sq B_{m0}$ or 
$\rau {P_m}\sg\sq B_{m1}$. 
Then $Q=\bigcap_mP_m$ is as required.
\epf

\parf{Proof of the main theorem}
\las{pmt}

\vyk{
If $m<\om$ and $\sg,\ta\in2^m$ then define    
$\zaf\sg\ta=
B\bez\ens{\phi(i)}{i<m\land \sg(i)\ne\ta(i)}$.%
\index{zzDsta@$\zaf\sg\ta$}%
}

Beginning the proof of Theorem~\ref{mt}, we let 
Borel maps ${\xe},{\xf}:\dn\to\dn$ witness the smoothness 
of the \eqr s resp.\ $\rE,\rF$, so that
$$
{x\rE y}\eqv {\xe}(x)={\xe}(y)
\qand
{x\rF y}\eqv {\xf}(x)={\xf}(y). 
$$
By Corollary \ref{bof}, we can assume that in fact 
${\xe},{\xf}$ \underline{are continuous}.

\ble
\lam{cap}
If\/ $P$ is a \pep, $U_0,U_1,\dots\sq \om$, 
and\/ ${\rE}$ is reduced to each\/ $U_k$ on\/ $P$, 
then\/ ${\rE}$ is reduced to\/ 
$U=\bigcap_kU_k$ on\/ $P$.
The same for\/ $\rF$.
\ele
\bpf
For just two sets, if $U=U_0\cap U_1$ and $x,y\in P$, 
$x\res U=y\res U$, then, using the product structure, 
find a point $z\in P$ with 
$z\res U_0=x\res U_0$ and $z\res U_1=y\res U_1$.
Then ${\xe}(x)={\xe}(z)={\xe}(y)$, hence $x\rE y$. 
The case of finitely many sets follows by  
induction. 
Therefore we can assume that 
$U_0\qs U_1\qs U_2\qs \ldots$ in the general case.
Let   
$x,y\in P$, $x\res U=y\res U$. 
There exist points $x_k\in P$  
satisfying $x_k\res U_k=x\res U_k$ and 
$x_k\res{(B\bez U_k)}=y\res {(B\bez U_k)}$. 
Then immediately ${\xe}(x_k)={\xe}(x)$, $\kaz k$.
On the other hand, clearly $x_k\to y$, hence, 
${\xe}(x_k)\to {\xe}(y)$ as ${\xe}$ is continuous. 
Thus ${\xe}(x)={\xe}(y)$, hence $x\rE y$.
\epf

\bpf[Theorem~\ref{mt}]
We argue in terms of Definition \ref{phi}. 
The plan is to define a sequence of \pep s  
as in Lemma \ref{fusm}, with some extra properties. 
Let $m<\om$. 
A \pep\ $R$ is \dd m\rit{good}, if 
(see definitions in Section~\ref{mud}):
\ben
\item[{\rm(1)$\rE$:}]
if $\sg\in2^m$ and $j=\phi(m)$ then either 
${\rE}$ is reduced to $\om\bez\ans j$ on $\rau R\sg$, 
or there is no \pep\ 
$R'\sq\rau R\sg$ on which
${\rE}$ is reduced to $\om\bez\ans j$;  

\item[{\rm(1)$\rF$:}]
the same for $\rF$;  

\item[{\rm(2)$\rE$:}]
if $\sg,\ta\in2^m$, then either 
(i) 
${\rE}$ is reduced on  
$\rau R\sg\cup\rau R\ta$ to 
$$
\zaf\sg\ta=
\om\bez\ens{\phi(i)}{i<m\land \sg(i)\ne\ta(i)}\,,
$$  
or 
(ii) 
$\imb {\xe}{\rau\zr\sg}\cap\imb {\xe}{\rau\zr\ta}=\pu$;\snos
{Given a function $h$ and $X\sq\dom h$, the set 
$\imb hX=\ens{h(x)}{x\in X}$ is the \dd himage of $X$.}  

\item[{\rm(2)$\rF$:}]
the same for $\rF$.
\een

\ble
\lam{ngel}
If\/ $m<\om$ and a \pep\/ $\zr$ 
is\/ \dd mgood, then there is an\/ 
\dd{m+1}good 
\pep\/ $\jq\nq{m+1}\zr$. 
\ele
\bpf[Lemma]
Consider a string $\sgi\in2^{m+1}$, 
and first define a    
\pep\ $\zq\in\mlt_B$, $\zq\nq{m+1}\zr$, 
satisfying (1)$\rE$ relatively to this string only.
Let $j=\phi(m+1)$. 
If there exists a \pep\ 
$\zr'\sq\rau\zr\sgi$ on which   
${\rE}$ is reduced to $\om\bez\ans j$, then  
let $\ju$ be such $R'$. 
If there is no such $R'$ then put 
$\ju=\rau\zr\sgi$. 
By Lemma \ref{nadm}, there is a \pep\ 
$\zq\nq{m+1}\zr$ such that $\rau\zq\sgi=\ju$. 
Thus the \pep\ $\zq$ satisfies (1)$\rE$   
with respect to $\sgi$.
Now take $\zq$ as the \lap{new} \pep\ $\zr$, 
consider another string $\sgi\in2^{m+1}$, and 
do the same as above. 
Consider all strings in $2^{m+1}$ 
consecutively, with the same procedure.
This ends with a \pep\ $\zq\nq{m+1}\zr$, 
satisfying (1)$\rE$ for all strings in $2^{m+1}$.

Now take care of (2)$\rE$. 
Let $\sgi,\tai\in2^{m+1}.$ 
Note that if $\sg'(m)=\ta'(m)$ then 
$\zaf{\sg'}{\ta'}=\zaf{\sg'\res m}{\ta'\res m}$, 
so that (2)$\rE$  relatively to $\sg',\ta'$ 
follows from (2)$\rE$  relatively to 
${\sg'\res m},{\ta'\res m}$. 
Thus it suffices to consider only pairs in $2^{m+1}$ 
of the form $\sg\we 0,\ta\we 1$, where $\sg,\ta\in2^m$.  
Consider one such a pair  
$\sgi=\sg\we 0$, $\tai=\ta\we 1$, 
and define a    
\pep\ $P\nq{m+1}\zq$, 
satisfying (2)$\rE$  relatively to this pair.

The sets $\cp=\zaf{\sgi}{\tai}$ 
and $U=\zaf{\sg}{\ta}$ satisfy 
$\cp=U\bez\ans{\eto}$, where $\eto=\phi(m)$, 
while the sets $\rau\zq\sgi$, $\rau\zq\tai$ 
satisfy 
$\rau\zq\sgi\res{\cp}=\rau\zq\tai\res{\cp}$. 

If ${\rE}$ is reduced to $U'$ 
on $Z'=\rau\zq\sgi\cup\rau\zq\tai$   
then (2)$\rE$(i) 
holds for $P=Q$ and the pair $\sgi,\tai$. 
Now suppose that ${\rE}$ is  
{\ubf not} reduced to $U'$ on $Z'$, 
so that there are points 
$x_0,y_0\in Z'$ with
$x_0\res U'=y_0\res U'$ and ${\xe}(x_0)\ne {\xe}(y_0)$, \ie, 
${\xe}(x_0)(k)=p\ne q={\xe}(y_0)(k)$ for some $k$ and  
$\ans{p,q}=\ans{0,1}$. 
As $\rau\zq\sgi\res{\cp}=\rau\zq\tai\res{\cp}$, 
we can \noo\ assume that $x_0\in \rau\zq\sgi$ 
but $y_0\in \rau\zq\tai$.      

As ${\xe}$ is continuous, there exist relatively clopen 
\pep s $X\sq \rau\zq\sgi$, $Y\sq \rau\zq\tai$, 
such that $x_0\in X$, $y_0\in Y$, 
${\xe}(x)(k)=p$ and ${\xe}(y)(k)=q$ for all $x\in X$, 
$y\in Y$. 
Define smaller \pep s $X'\sq X$ and $Y'\sq Y$ so 
that $X'(j)=Y'(j)=X(j)\cap Y(j)$ for all $j\in U'$ 
but $X'(j)=X(j)$, $Y'(j)=Y(j)$ 
for all $j\in\om\bez U'$. 
Note that still $x_0\in X'$, $y_0\in Y'$, and now 
$X'\res U'=Y'\res U'$.

This allows to apply Lemma~\ref{nadm}\ref{nadm4}, 
getting a \pep\ $P\nq{m+1} Q$ such that 
$\rau P\sgi=X'$ and $\rau P\tai=Y'$.
Then $\imb {\xe}{\rau P\sgi}\cap\imb {\xe}{\rau P\tai}=\pu$
by construction, therefore (2)$\rE$(ii) 
holds for $P$ and the pair of $\sgi,\tai$. 

To conclude, we get a \pep\ $P\nq{m+1} Q$ such that 
(2)$\rE$ holds for $P$ and the pair of $\sgi,\tai$
in both cases. 

Consider all pairs of strings in $2^{m+1}$ 
consecutively.
This yields a \pep\ $R\nq{m+1}Q$, 
satisfying (2)$\rE$ for all $\sgi,\tai\in2^{m+1}$ 
(and still satisfying (1)$\rE$).

Then repeat the same procedure for $\rF$.
\epF{Lemma}

Come back to the proof of the theorem. 
Lemma \ref{ngel} yields  
an infinite sequence 
$\dots\leq_3\zq_2\leq_2\zq_1\leq_1\zq_0$
of \pep s $\zq_m$, such that each $\zq_m$ 
is a \dd{m}good.
The limit \pep\ 
$P=\bigcup_m\zq_m\in\mlt_B$ satisfies  
$P\nq{m+1}\zq_m$ for all $m$ by Lemma \ref{fusm}. 
Therefore $P$ is \dd{m}good for every $m$, hence 
we can freely use (1)$\rE,\rF$ and (2)$\rE,\rF$ 
for $P$ in the following final argument.\vom

{\ubf Case 1\/}: 
if $m<\om$, $\sg,\ta\in2^m$, and 
$\imb {\xe}{\rau P\sg}\cap\imb {\xe}{\rau P\ta}=\pu$, 
then 
$\imb {\xf}{\rau P\sg}\cap\imb {\xf}{\rau P\ta}=\pu$.
Prove that ${\rF}\sq{\rE}$ on $P$ 
in this case, as required by the 
``either'' option of Theorem~\ref{mt}.
Assume that $x,y\in P$ and $x\rE y$ {\ubf fails}, 
that is, ${\xe}(x)\ne {\xe}(y)$; 
show that ${\xf}(x)\ne {\xf}(y)$.
Pick $a,b\in\dn$ satisfying 
$\ans x=\bigcap_m \rau P{a\res m} $ and 
$\ans y=\bigcap_m \rau P{b\res m} $. 
As $x\ne y$, we have 
$\imb {\xe}{\rau\zq{a\res m}}\cap
\imb {\xe}{\rau\zq{b\res m}}=\pu$
for some $m$ 
by the continuity and compactness. 
Then by the Case 1 assumption, 
$\imb {\xf}{\rau P{a\res m}}\cap
\imb {\xf}{\rau P{b\res m}}=\pu$ holds, 
hence ${\xf}(x)\ne {\xf}(y)$, and $x\rF y$ {\ubf fails}.\vom

{\ubf Case 2\/} = not Case 1. 
Then, by (2)$\rF$, 
there is a pair of strings
$\sgi=\sg\we i, \:\tai=\ta\we k\in2^{m+1}$, $m<\om$, 
such that 
$\imb {\xe}{\rau P\sgi}\cap\imb {\xe}{\rau P\tai}=\pu$, 
but ${\rF}$ is reduced to $U'=\zaf\sgi\tai$ 
on $Z'=\rau P\sgi\cup\rau P\tai$.
Assume that $m$ is the least possible 
witness of this case.
We are going to prove that the \pep\ $\rau P\sg$ 
satisfies the ``or'' option of Theorem~\ref{mt}, 
with the number $\eto=\phi(m)$, 
that is, 
(*) 
${\rF}$ is reduced to $\om\bez\ans\eto$
on $\rau P\sg$, 
and 
(**) 
${\rE}$ captures $\eto$  
on $\rau P\sg$. 

\ble
\lam{nge2}
The relation\/ ${\rE}$ is$:$\vom 

{\rm(A)} 
reduced to\/ $U=\zaf\sg\ta$ on the set\/  
$Z=\rau P\sg\cup\rau P\ta$, 
\vom 

{\rm(B)} 
{\ubf not} reduced to\/ $U'=\zaf\sgi\tai$ on\/  
$Z'=\rau P\sgi\cup\rau P\tai$,\vom 

{\rm(C)} 
{\ubf not} reduced to 
$\om\bez\ans\eto$ 
on any \pep\ $\ju\sq\rau{P}\sg$.\vom

\noi
In addition, {\rm(D)} $U\ne U'$, hence\/ 
$\eto\in U$ and\/ $U'=U\bez\ans\eto$.
\ele
\bpf
(A) 
Otherwise  
$\imb {\xe}{\rau P\sg}\cap\imb {\xe}{\rau P\ta}=\pu$ 
by (2)$\rE$, 
hence ${\rF}$ is {\ubf not} reduced to $U$  on 
%
$\rau P\sg\cup\rau P\ta$ by the choice of $m$, 
thus 
$\imb {\xf}{\rau P\sg}\cap\imb {\xf}{\rau P\ta}=\pu$ 
by (2)$\rF$, then 
$\imb {\xf}{\rau P\sgi}\cap\imb {\xf}{\rau P\tai}=\pu$, 
which contradicts to the fact that 
${\rF}$ is reduced to $U'$ on 
$\rau P\sgi\cup\rau P\tai$, 
as  
$\rau P\sgi\res U'=\rau P\tai\res U'$ 
by Lemma~\ref{nadm}.

(B) 
The otherwise assumption contradicts to 
$\imb {\xe}{\rau P\sgi}\cap\imb {\xe}{\rau P\tai}=\pu$. 

(D) follows from (A) and (B).

(C) 
Otherwise 
${\rE}$ is reduced to $\om\bez\ans\eto$ 
on $ \rau{P}\sg $ by (1)$\rE$. 
Then ${\rE}$ is reduced to $U'$ on $ \rau{P}\sg$ 
by Lemma \ref{cap} since $U'=U\bez\ans\eto$ by (D).
%
It follows that ${\rE}$ is reduced to $U'$     
on $Z$,\snos 
{Let $x,y\in Z= \rau{P}\sg \cup \rau{P}\ta$ 
and $x\res U'=y\res U'$. 
As $\rau{P}\sg\res U=\rau{P}\ta\res U$ 
by Lemma \ref{nadm}, 
there are $x',y'\in \rau{P}\sg $ with  
$x\res U=x'\res U$ and $y\res U=y'\res U$.
We have $x\rE x'$ and $y\rE y'$ 
by (A), and $x'\rE y'$ since 
${\rE}$ is reduced to $U'$ on $\rau{P}\sg $.
We conclude that $x\rE y$.
}
hence on $Z'\sq Z$ 
as well.
But this contradicts to (B). 
\epf


Now, as $U'=U\bez\ans\eto\sq\om\bez\ans\eto$, 
the \pep\ $\rau P\sgi$
witnesses that ${\rF}$ is reduced to $\om\bez\ans\eto$
%
on $\rau P\sg$ by (1)$\rF$.
Thus we have (*).

To prove (**), 
let $x,y\in  \rau{P}\sg $ and $x\rE y$; 
prove that $x(\eto)=y(\eto)$. 
Indeed we have    
$\ans x=\bigcap_n \rau{P}{a\res n} $ and  
$\ans y=\bigcap_n \rau{P}{b\res n} $, 
where $a,b\in\dn$, $\sg\su a$, $\sg\su b$. 
Let $\zaf ab=\bigcap_n\zaf{a\res n}{b\res n}$. 
Then $x\res\zaf ab=y\res\zaf ab$, since 
$\rau{P}{a\res n}\res{\zaf{a\res n}{b\res n}}=
\rau{P}{b\res n}\res{\zaf{a\res n}{b\res n}}$
for all $n$. 
Thus it suffices to check that 
$\eto\in\zaf{a\res n}{b\res n}$ for all $n$.

Suppose towards the contrary that 
$\eto=\phi(m)\nin\zaf{a\res n}{b\res n}$ 
for some $n$. 
Note that $n>m$ because $a\res m=b\res m=\sg$. 
However ${\rE}$ is reduced to 
$\zaf{a\res n}{b\res n}$
on $ \rau{P}{a\res n}$
by (2)$\rE$, since $x\rE y$. 
Yet we have $\eto\nin\zaf{a\res n}{b\res n}$, 
therefore, $\zaf{a\res n}{b\res n}\sq \om\bez\ans\eto$. 
It follows that ${\rE}$ is reduced to $\om\bez\ans\eto$
on $ \rau{P}{a\res n}$. 
But this contradicts to Lemma~\ref{nge2}(C) with 
$\ju=\rau{P}{a\res n}$.  

To conclude Case 2, we have checked (*) and (**). 
\epF{Theorem~\ref{mt}}

\bibliographystyle{plain}

\end{document}